\documentclass[10pt]{article}
\usepackage[usenames]{color}
\definecolor{Red}{rgb}{1,0,0}
 \usepackage[parfill]{parskip}   
\usepackage{graphicx}
\usepackage{amssymb}
\usepackage{epstopdf}
 \usepackage{times}
\usepackage{amsthm}
\usepackage{amsmath}
\usepackage{amscd}
\usepackage{makeidx}
\usepackage{graphicx}
 \newtheorem{definition}{Definition}[section]
\newtheorem{proposition}{Proposition}[section]

\newtheorem{lemma}[proposition]{Lemma}
\newtheorem{theorem}[proposition]{Theorem}

\newtheorem{remark}{Remark}[section]
\newtheorem{assumption}{Assumption}[section]

\begin{document}
  \title {On the universal principally polarized abelian variety of dimension  4
 }
\author{  Alessandro Verra \\ Dipartimento di Matematica, Universit\'a Roma Tre }  
\date {}
    \maketitle
 \section{Introduction}
 
 Let   ${\cal A}_g$ be the coarse moduli space for principally polarized abelian varieties of dimension $g$ over an algebraically closed field $k$, $char \ k \neq 2,3$.  In this paper we deal with  the universal family
$$
f: {\cal {X}}_g \to {\cal A}_{g},
$$
$\textcolor{Red}{{\cal T}_g} \subset \textcolor{Red}{{\cal X}_g}
$
defined by the following property: for $u \in U$ let $A = f^*u$  and let $\Theta = A \cdot \textcolor{Red}{{\cal T}_g}$. Then $\Theta$ is a principal polarization on $A$ and $u$ is the moduli point of the
p.p.a.v. $(A, \Theta)$. We will say that $\textcolor{Red}{{\cal T}_g}$ is the \it universal theta divisor on \textcolor{Red}{\textcolor{Red}{$\textcolor{Red}{{\cal A}_g}$}}. \rm  \\
In some sense the  study of  \textcolor{Red}{$\textcolor{Red}{{\cal X}_g}$}  parallels  that of   the universal  Jacobian ${\cal Y}_g$  over the moduli space $\textcolor{Red}{{\cal M}}_g$ of curves of genus $g$.  For low $g$ some natural questions are still open: for instance one would like to know the Kodaira dimension of \textcolor{Red}{$\textcolor{Red}{{\cal X}_g}$} and ${\cal Y}_g$ for any $g$. Let us discuss a related problem:  \par \noindent \it  
For which values of $g$ is \textcolor{Red}{$\textcolor{Red}{{\cal X}_g}$} uniruled, rationally connected, unirational or even rational? \rm \par \noindent  Since $f$ is dominant none of these  properties appears if  the Kodaira dimension of \textcolor{Red}{\textcolor{Red}{$\textcolor{Red}{{\cal A}_g}$}} is $\geq 0$,  in particular if  \textcolor{Red}{\textcolor{Red}{$\textcolor{Red}{{\cal A}_g}$}} is of general type. This is actually the case for $g \geq 7$ due to the  results of Freitag, Mumford and Tai ([F], [M1], [T]).   So we are confined to $g \leq 6$. \par \noindent 
If $g \leq 3$ then \textcolor{Red}{$\textcolor{Red}{{\cal X}_g}$} is birational to ${\cal Y}_g$. The rationality of $\textcolor{Red}{{\cal X}_1}$ is well  known,
while the results about the rationality of  moduli of pointed curves, [CF], should also imply the rationality of $\textcolor{Red}{{\cal X}_2}$ and $\textcolor{Red}{{\cal X}_3}$.
\par \noindent
The genus 4, 5 and 6 cases are open. Notice also that \textcolor{Red}{${\cal A}_4$} and $\textcolor{Red}{{\cal A}_5}$ are known to be unirational, (see the survey [G] and [C], [D], [L], [ILS], [V] ), while the Kodaira dimension of $\textcolor{Red}{{\cal A}_6}$ is unsettled. In this note we contribute to the genus 4 case of the previous problem  proving the following:  \par \noindent
 {\bf Theorem} \  \it (1) The universal principally polarized abelian variety over \textcolor{Red}{${\cal A}_4$}  is unirational. \\
 (2) The universal theta divisor over \textcolor{Red}{${\cal A}_4$} is unirational.
\endtheorem \rm  \par \noindent
The proof relies on the theory of Prym varieties,  on the beautiful geometry related to  \'etale double covers of curves of low genus and on K3 surfaces endowed with a Nikulin involution. We conclude this introduction with a brief  summary of it. \par \noindent  
Due to the Prym construction the universal family ${\cal X}_4$ is dominated by the moduli space of pairs
$$
(\pi, d)
$$
where $\pi: \tilde C \to C$ is a non split  \'etale double cover of a smooth, irreducible curve of genus 5 and $d$ is an isolated, effective divisor on $\tilde C$ such that $\pi_* d \in \mid \omega_C \mid$.   We prefer the moduli of triples 
$$
(\pi, d, L)
$$
where $L = {\cal O}_{\tilde C}(l) \in Pic^8(\tilde C)$ satisfies $h^0( L) = 3$ and $\pi_*l \in \mid \omega_C \mid$.
By  the Brill-Noether theory for special curves  $\tilde C$ as above a 1-dimensional family of line bundles $ L$ does exist on $\tilde C$. Therefore this new moduli space dominates ${\cal X}_4$ via the forgetful map  $( \pi, d, L)$ $\to$ $(\pi, d)$. \par \noindent The additional line bundle $ L$ is geometrically useful:  indeed let
$
i: \tilde C \to \tilde C
$
be the fixed point free involution induced by $\pi$, it turns out that $ L \otimes i^*  L \cong \omega_{\tilde C}$ and that the Petri map
$$
\mu: H^0( L) \otimes H^0(i^*L) \to H^0(\omega_{\tilde C})
$$
has corank one for a general $\tilde C$ as above. This is an important point: with a little extra work one can deduce from this property that the linear system $\mid Im \ \mu \mid$ defines an embedding
$$
d \subset  \tilde C   \subset \mathbf P^2 \times \mathbf P^2  \subset \mathbf P^8,
$$
where the linear span of $\tilde C$ is a hyperplane in $\mathbf P^8$. Let $\iota$ be the  projective involution exchanging  the factors of $\mathbf P^2 \times \mathbf P^2$, in addition we can assume that $i = \iota / \tilde C$.  In sections 3, 4 and 5 we use this embedding to construct a rational parameter space for the family of the $0$-cycles $d$ as above. Essentially we construct a rational  family of triples $ (o, \tilde S, d)$ with the following properties: \par \noindent
(1) $o = (y, B)$ where $B \subset \mathbf P^2 \times \mathbf P^2$ is a smooth conic and $y = (a, b_1,b_2)$ is a triple such that:  $(b_1, b_2) \in B \times B$, $a$ is a set of 6 independent points in $\mathbf P^2 \times \mathbf P^2$, $q(a)$ is a set of 6 coplanar points where $q: \mathbf P^8 \to \mathbf P^5$ is the projection from the 2-dimensional  projectivized eigenspace of $\iota$.  \par \noindent
(2)  $\tilde S$ is a smooth K3 surface  which is a  complete intersection of a quadratic and a hyperplane section of $\mathbf P^2 \times \mathbf P^2$. Moreover $\iota$ restricts to a Nikulin involution on $\tilde S$ and $a \ \cup \ B \subset \tilde S$.  \par \noindent
(3)  \  $a \ \cup \ b_1 \ \cup \ b_2$ is contained in a unique, smooth hyperplane section $\tilde F$ of $\tilde S$. Moreover $\iota$ restricts to a fixed-point-free involution on $\tilde F$ and finally 
$$
d \in \mid {\cal O}_{\tilde F}(a + b_1 + \iota(b_2) ) \mid.
$$
We will say that $(o, \tilde S, d)$   is a \it  marking $0$-cycle of type 1. \rm  It turns out that the latter linear system is a pencil of divisors of degree 8 and that the family
$\mathbb D_1$
of marking 0-cycle of type 1  is rational. \par \noindent  Given an element  $ (o, \tilde S, d) \in \mathbb D_1$   we have from it  the pencil $\Lambda$ of all curves 
$$ \tilde C \in \mid {\cal O}_{\tilde S}(\tilde F + B + \iota^*B) \mid $$ passing through $d$ and having $\iota^*$-invariant equation. As we will see a general
$
\tilde C \in \Lambda
$
is a smooth, irreducible curve of genus 9. In addition $\tilde C$ is endowed with the fixed-point-free involution $\iota / \tilde C$ and it is marked by the divisor $d$. Let $C = \tilde C / <\iota>$ and let $\pi: \tilde C \to C$ be the quotient map,  it turns out that $\pi_*d \in \mid \omega_C \mid$ and that $h^0({\cal O}_{\tilde C}(d)) = 1$. Hence ${\cal O}_{\tilde C}(d)$ defines a point of the Prym variety $P(\pi)$ of $\pi$. Since $P(\pi)$ is a 4-dimensional p.p.a.v. the pair $(\pi,d)$ defines a point of the universal family ${\cal X}_4$. Finally we consider the family $\tilde {\mathbb C}_1$ of all 4-tuples $(o, \tilde S, d, \tilde C)$ and the map
$$
\phi_1: \tilde {\mathbb C}_1 \to {\cal X}_4,
$$
sending  $(o, \tilde S, d, \tilde C)$ to the moduli point of the pair $(P(\pi),d)$. $\tilde {\mathbb C}_1$ is a $\mathbf P^1$-bundle over $\mathbb D_1$, hence it is rational. We show in section 6 that $\phi_1$ is dominant, so that ${\cal X}_4$ is unirational.  \par \noindent
The proof of the unirationality of the universal theta divisor over \textcolor{Red}{${\cal A}_4$} is identical: in this case we use the family of \it marking 0-cycles of type 2. \rm These are triples $(o, \tilde S,d)$ defined by the same conditions (1), (2), (3):  the only difference is that
$
d \in \mid {\cal O}_{\tilde F}(a + b_1 + b_2) \mid.
$
We omit further details. \medskip \par
\bf Aknowledgements \rm Prym varieties are nowadays a subject having a long recent history, starting perhaps in the early seventies of the last century. One of the main actors of this  history is certainly Roy
Smith,  whose 65-th birthday is celebrated in this volume. So it is a pleasure and a honour, for  the author of  this paper, to contribute and to wish all the best to Roy. \par \noindent
\remark \rm As the referee pointed out,  the results of this paper imply that the Prym moduli space \textcolor{Red}{$\textcolor{Red}{{\cal A}_5}$} is unirational. Indeed  \textcolor{Red}{$\textcolor{Red}{{\cal A}_5}$} is the moduli space of  \'etale double covers $\pi: \tilde C \to C$ as above and it turns out that the map
$
\rho: \tilde {\mathbb C}_1 \to \textcolor{Red}{{\cal A}_5},
$
sending $(o, \tilde S, d, \tilde C)$ to the moduli point of $\pi$, is dominant. \\ About this one has however to mention that an independent and simple proof of the unirationality of \textcolor{Red}{$\textcolor{Red}{{\cal A}_5}$} was recently worked out by Marco Lo Giudice in his Tesi di Dottorato, ([L]). An amplified and revised version of such a proof  is contained in the preprint   \it The moduli space of \'etale double covers of genus 5 curves is unirational, \rm by E. Izadi, M. Lo Giudice and G. Sankaran, ([ILS]). \endremark \rm \par \noindent
At the end of this introduction we wish to thank the referee: for his previous  remark on \textcolor{Red}{$\textcolor{Red}{{\cal A}_5}$}, for some useful observations and finally for his careful and patient reading of the paper.

\endremark
 
\section {Preliminaries on Pryms and notations}
Usually we will denote as
$$
\pi: \tilde C \to C
$$
a non split \'etale double cover of a smooth, irreducible curve $C$ of genus $g$ and by
$$
i: \tilde C \to \tilde C
$$
the fixed-point-free involution exchanging the sheets of $\pi$. To give
$\pi$ is equivalent to giving a non trivial order two element $\eta \in Pic^0(C)$.  The Prym variety associated to $\pi$  will be denoted by
$$
P(\pi),
$$
as is well known $P(\pi)$ is a principally polarized abelian variety of dimension $g-1$. Let us
recall a useful construction of $P(\pi)$ described by Mumford in [M2]: consider the Norm map
$$
Nm: Pic^{2g-2}(\tilde C) \to Pic^{2g-2}(C),
$$
sending ${\cal O}_{\tilde C}(\sum x_i)$ to ${\cal O}_C(\sum \pi(x_i))$. Each fibre
of $Nm$ is the disjoint union of 2 copies of $P(\pi)$. For the point $\omega_C \in Pic^{2g-2}(C)$
it turns out that $Nm^{-1}(\omega_C) = P^+ \cup P^-$, where
$$
P^+ = \lbrace M \in Pic^{2g-2} (\tilde C) \ / \ Nm( M) = \omega_C \ , \ h^0(M) \ is \ even \rbrace
$$
and
$$
P^- = \lbrace M \in Pic^{2g-2} (\tilde C) \ / \ Nm( M) = \omega_C \ , \ h^0(M) \ is \ odd \rbrace.
$$
For $M$ general in $P^-$ one has $h^0(M) = 1$ and one has $h^0(M) = 0$ for $M$ general in $P^+$. Let
$$
\Xi = \lbrace M \in P^+ \ / \ h^0(M) \geq 2 \rbrace,
$$
then $\Xi$ is a principal polarization on $P^+$: by definition $P(\pi)$ is the pair $(P^+, \Xi)$. We will 
use the following well known property, sometimes called \it Parity Lemma, \rm (cfr. [B] prop. 3.4 and [M3]).
\lemma Let $x \in \tilde C$ then $M \in P^-$ if and only if $M(x-i(x)) \in P^+$.
\endlemma
\proof  Let $n \in Div^{2g-2}(\tilde C)$ be an effective divisor such that $N := {\cal O}_{\tilde C}(n)\in P^-$. If $N$ is general in $P^-$ then $x$, $i(x)$ are not in $Supp \ n$ and $h^0(N(x)) = 1$. Therefore  $h^0(N(x-i(x)) = 0$ and  the translation $M \to M(x-i(x))$ induces a biregular isomorphism between $P^-$ and $P^+$. This implies  the lemma. \endproof \rm \par \noindent
 The  Brill-Noether theory  is  known for $\tilde C$, though $\tilde C$ is not a general curve ([W]).  Let 
$$
\mu: H^0(M) \otimes H^0(i^*M) \to H^0(\omega_{\tilde C})
$$
be the Petri map, for a given $M \in P^+ \cup P^-$. Then $i^*$ acts as an involution on the above vector spaces and $\mu$ preserves eigenspaces. By definition the Prym-Petri map
$$
\mu_-: [ H^0(M) \otimes H^0(i^*M)]^- \to H^0(\omega_{\tilde C})^-
$$
is the induced map between the corresponding $-1$ eigenspaces.  As is well known $\mu^-$ is injective for a general $\tilde C$ as above and every $M$. Moreover consider the Prym-Brill-Noether scheme
$$
W^r(\pi) := \lbrace N \in P^+ \ (P^-)  \ / \ h^0(N) \geq r+1 \ , \ h^0(N) = r + 1  \ mod \ 2 \rbrace.
$$ 
For any $\tilde C$ as above $Coker \ \mu^-$ is the tangent space to  $W_r(\pi)$ at $M$. For  general $\tilde C$ and $M$ the rank of  $\mu^-$ is $\binom {r+1}2$ so that $dim \ Coker \ \mu^- = g - 1 -  \binom {r+1}2$.  \par \noindent Finally we recall that the Prym moduli space
$
\textcolor{Red}{{\cal R}_g}
$
is by definition  the moduli space of the \'etale double covering $\pi$. \textcolor{Red}{$\textcolor{Red}{{\cal R}_g}$} is irreducible. Let
$
p_g: \textcolor{Red}{{\cal R}_g} \to \textcolor{Red}{{\cal A}_{g-1}}
$
be the Prym map sending the moduli point of $\pi$ to the moduli point of $P(\pi)$. Then
$p_g$ is dominant for $g \leq 6$. In particular its  degree is 27 for $g = 6$, (see [DS]). \par \noindent 
\underline {\it Some frequent notations:} \rm  \vfill \eject  \noindent
$\circ$  $<S>$ will denote the linear span of a subset $S \subset \mathbf P^n$. \\
$\circ$  $V^+$ and $V^-$ will denote the +1 and -1 eigenspaces of a fixed involution $i$ of $V$. \\
$\circ$ For $\Lambda = \mathbf PV$ we will use the notations $\Lambda^+$ and $\Lambda^-$ for the projectivized eigenspaces of $i$.\\
$\circ$ For the ideal sheaf of $Y$ in $X$ we will use the notation ${\cal I}_{Y/X}$.
 \section {Conics and some useful  0-cycles in $\mathbf P^2 \times \mathbf P^2$}
First, we fix some further  notations and conventions: we begin by fixing an \it auxiliary projective plane 
$
\Pi.
$
\rm Let  $V := H^0({\cal O}_{\Pi}(1))$ and $\mathbf P^2 := \mathbf PV$ we will consider the Segre embedding
$$
\mathbf P^2 \times \mathbf P^2 \subset  \mathbf P^8 := \mathbf P (V \otimes V).
$$
The natural projection of $\mathbf P^2 \times \mathbf P^2$ onto the $i$-th factor will be  denoted as
$$
\pi_i:  \mathbf P^2 \times \mathbf P^2 \to \mathbf P^2, \ \ i = 1,2.
$$
From the decomposition $V \otimes V = \wedge^2  V \oplus Sym^2 V$ we obtain the subspaces
$$
\mathbf P^- := \mathbf P \wedge^2 V = \Pi  \  \rm and \it \ \mathbf P^+ := \mathbf P Sym^2 \ V = \mid {\cal O}_{\Pi}(2) \mid
$$
which are the projectivized eigenspaces of the involution $v_1 \otimes v_2 \to v_2 \otimes v_1$. We denote by
$$
\iota: \mathbf P^8 \to \mathbf P^8
$$
the induced projective involution, while  the natural linear projections onto $\mathbf P^-$ and $\mathbf P^+$ will be 
$$
p: \mathbf P^8 \to \mathbf P^- \ \rm and \  \it q: \mathbf P^8 \to \mathbf P^+.
$$
It is clear that $\mathbf P^2 \times \mathbf P^2$ is the space of ordered pairs $(l,l')$ of lines of $\Pi$ and that $i(l,l') = (l',l)$.  It turns out that: \ (1) $p(l,l')$ is the point $x = l \cap l'$ of $\Pi$, \ (2) $q(l,l')$ is the singular conic $l+l'$ of $\Pi$. \ Notice that
$$
\Delta := \mathbf P^2 \times \mathbf P^2 \cap \mathbf P^+
$$
is the diagonal of $\mathbf P^2 \times \mathbf P^2$. $\Delta$ is  embedded as a Veronese surface
in the 5-dimensional space $\mathbf P^+$. The variety of the bisecant lines to $\Delta$ is a well known cubic hypersurface, we denote it by 
$$
\Sigma.
$$
Finally we mention another well known fact i.e.  that the following diagram is commutative:
$$
\CD
{\mathbf P^-} @<p<< {\mathbf P^8} @>q>> {\mathbf P^+}  \\
{\cup} @. {\cup} @. {\cup} \\
{\mathbf P^-} @<{p/ \mathbf P^2 \times \mathbf P^2}<< {\mathbf P^2 \times \mathbf P^2} @>{q/ \mathbf P^2 \times \mathbf P^2} >> {\Sigma} \\
 \endCD
$$
In particular $q/\mathbf P^2 \times \mathbf P^2$ is a finite 2:1 cover of $\Sigma$ branched along $\Delta$ and $\Delta = Sing \ \Sigma$. \par \noindent
Now we begin our constructions by defining a suitable family of $0$-dimensional subschemes of length 6
in $\mathbf P^2 \times \mathbf P^2$. In the Grassmannian of 4-spaces of $\mathbf P^8$ let $\sigma $ be the  Schubert cycle parametrizing all 4-spaces $\Lambda$ such that $dim \Lambda \cap \mathbf P^- \geq 1$. Then $\sigma$ is rational and contains the open set 
$$
A := \lbrace \Lambda \in \sigma \ /  \ \text { $\Lambda$ is transversal to $\mathbf P^2 \times \mathbf P^2$,  \  $\Lambda \  \cap \mathbf P^-$ is a line,  \  $\Lambda \cap \iota(\Lambda) \cap \mathbf P^2 \times \mathbf P^2 = \emptyset$} \rbrace.
$$
For each $ \Lambda \in A$ the scheme $Z = \Lambda \cdot \mathbf P^2 \times \mathbf P^2$ has  length 6 and spans $\Lambda$. Moreover  $\Lambda = < Z - z >$ for each $z \in Z$: this just follows because $\mathbf P^2 \times \mathbf P^2$ is a 4-fold of minimal degree 6 in $\mathbf P^8$. The definition of $A$ also implies that $$q/Z: Z \to q(Z)$$ is bijective:  if not we would have $q(x) = q(y)$ for two distinct points $x,y \in Z$  and hence  $y = \iota(x)$. This implies $x, y \in \Lambda \cap \iota(\Lambda) \cap \ \mathbf P^2 \times \mathbf P^2$,  which is excluded.  From now on we identify $A$ to the family of schemes $Z = \Lambda \cdot \mathbf P^2 \times \mathbf P^2$, where $\Lambda \in A$.  Over $A$ we consider the universal family: \par \noindent
\definition $\mathbb A_0:= \lbrace (Z,z) \in  A \times \mathbf P^2 \times \mathbf P^2 \ / \ z \in Z \rbrace$.
\par \noindent
 \proposition $\mathbb A_0$ is rational. \endproposition
\proof Consider the map $\phi: \mathbb A_0 \to \mathbf P^{-*} \times \mathbf P^2 \times \mathbf P^2$ defined as follows. Let $(Z,z) \in \mathbb A_0$, then $z$ is a point of the scheme $Z \in U$ and $l_Z := <Z> \cap \ \mathbf P^-$ is a line: by definition $\phi(Z,z) := (l_Z,z)$. Note  that the fibre of $\phi$ at $(l_Z,z)$ is open in the Grassmann variety $G(1,5)$ parametrizing all the 4-spaces containing the  plane $<l_Z,z>$. It is standard to deduce from this that $\mathbb A_0$ is open in a $G(1,5)$-bundle over $\mathbf P^{-*} \times \mathbf P^2 \times \mathbf P^2$: we omit the details. Then $\mathbb A_0$ is rational. \endproof \rm \par \noindent
Let $(Z,z) \in \mathbb A_0$: in $Z$ we can replace $z$ by $\iota(z)$, obtaining a second scheme we denote as  
$$
Z'.
$$
\definition $\mathbb A := \lbrace (Z', \iota(z)), \ (Z,z) \in \mathbb A_0 \rbrace$.
\enddefinition \rm \par \noindent
It is clear from the definition that  $\mathbb A$ is biregular to $\mathbb A_0$, so that $\mathbb A$ is rational too. 
\proposition Let $(Z,z) \in {\mathbb A}_0$. Then the following properties hold:  \par \noindent
(1)  $<Z'>$ is a 5-dimensional space containing $\mathbf P^-$, \par \noindent
(2) $<Z>$ is a 4-dimensional space and $<Z> \ \cap \  \mathbf P^-$ is a line,\par \noindent
(3) $<q_*Z> = <q_*Z>$, moreover this is a plane in $\mathbf P^+$.
\endproposition \rm \par \noindent
\proof (1)  Recall that $<Z>$ is a 4-space intersecting $\mathbf P^-$ along a line and not containing $\iota(z)$. Since $<Z> = <Z-z>$ it follows  that $<Z'>$ is the 5-space spanned by $Z$  and $\iota(z)$. 
In particular $<Z'>$ contains the line $ l =<z \ \iota(z)>$ and $l$ intersects  $\mathbf P^-$ in a point not in $<Z> \cap \ \ \mathbf P^-$. Hence it follows $\mathbf P^- \subset \ <Z'>$. (2) and (3) are obvious consequences of  the definitions and of (1). \endproof \rm \par \noindent
For simplicity  we will use the \it notation \rm $a$ both for an element of $\mathbb A \cup \mathbb A_0$ and for its corresponding scheme $Z$ or $Z'$,   omitting to indicate the distinguished point of such a scheme unless it is necessary.    \par \noindent
 \definition Let $a \in \mathbb A$, then  $\mathbf P^5_a := \  < a >$ and \ $\mathbf P^{2+}_a \  := \mathbf P^5_a \cap \mathbf P^+$. \enddefinition \rm \par \noindent
\remark \rm  $\mathbf P^{2+}_a$ is  the plane spanned by  $q_*a$.  By proposition 3.2 (1) $\mathbf P^5_a$ contains $\mathbf P^-$. Then $\mathbf P^5_a$ is spanned by $\mathbf P^{2+}_a$ and $\mathbf P^-$ and hence it  is $\iota$-invariant: the set of fixed points of $\iota / \mathbf P^5_a$ is $\mathbf P^- \cup \mathbf P^{2+}_a$. \endremark \rm \par \noindent

Let $P$ be a general $5$-space through $\mathbf P^-$:  it is clear that $P$ is transversal to $\mathbf P^2 \times \mathbf P^2$ and that $P = \mathbf P_a$ for some $a \in \mathbb A$. This implies that $\mathbf P^5_a$ is transversal to $\mathbf P^2 \times \mathbf P^2$ for a general $a \in \mathbb A$. \par \noindent
\assumption We always assume that $a$ is general in the above sense. 
\endassumption  \rm \par \noindent
Now we consider in $\mathbb A \times \mathbf P^2 \times \mathbf P^2$ the open set $U$ of pairs
$(a,b_1)$ such that the linear span $<a \cup b_1>$ is 6-dimensional and transversal to $\mathbf P^2 \times \mathbf P^2$. This implies that for each $u = (a,b_1) \in U$  
$$
Y_u := <a \cup b_1> \cap \  \mathbf P^2 \times \mathbf P^2
$$
is a smooth sextic Del Pezzo surface. Then we consider over $U$  the universal family
$$
\mathbb Y= \lbrace (u,b_2) \in U \times \mathbf P^2 \times \mathbf P^2 \ / \ b_2 \in Y_u \rbrace.
$$
\proposition $\mathbb Y$ is rational. \endproposition
\proof Let $p: \mathbb Y \to U$ be the natural projection and let $\alpha: U \to \mathbb Y$ be the section
which is so defined: if $u = (a,b_1)$,  $\alpha(u)$ is the distinguished point of $a$. For each  $u = (a,b_1)$ consider the map
$$
\phi_u: Y_u \to \mathbf P^2
$$
defined by the linear system of hyperplane sections of $Y_u$ which are singular at $b_2$ and contain
the point $\alpha(u)$. It is easy to see that this map is birational and to deduce that there exists a
birational map $\phi: \mathbb Y \to U \times \mathbf P^2$ such that $\phi/ Y_u = \phi_u$. Since $\mathbb A$
is rational, this implies that $\mathbb Y$ is rational. \endproof \rm \par \noindent
With some abuse we will still denote by $\mathbb Y$  its open set defined by the following condition (*): \par \noindent
(*) let $(u,b_2) \in \mathbb Y$ with $u = (a,b_1)$ then: \par \noindent
(1) $b_1, b_2$ are distinct points not in $\mathbf P^5_a \cup \Delta$ and such that $b_2 \neq \iota (b_1)$, \par \noindent
(2) their projections $\pi_i(b_1)$ and $\pi_i(b_2)$ are distinct points in $\mathbf P^2$, for $i = 1,2$.\par \noindent
(3) Any 7 points of $ a \cup \lbrace b_1, \iota(b_1), b_2, \iota(b_2) \rbrace$, no two exchanged by $\iota$, are linearly independent.
\par \noindent
\remark \rm For completeness we show that (*) is satisfied on a non empty open set of $\mathbb Y$: \\
Let $u = (a, b_1) \in \mathbb Y_u$, it is clear that $(u,b_2)$ satisfies (1) and (2) if $b_2$ is  general in $\mathbb Y_u$. It follows from (1) that any $x \in \beta := \lbrace b_1, b_2, \iota(b_1), \iota(b_2) \rbrace$ is not  in $\mathbf P^5_a$, so that the seven points of $a \cup x$ are linearly independent. Let $(x, y)$ be
one of the following pairs: $(b_1, b_2)$, $(b_1, \iota(b_2))$, $(\iota (b_1), b_2)$, $(\iota(b_1), \iota(b_2))$. Since $(b_1, b_2)$ is general in $\mathbb Y_u \times \mathbb Y_u$, the same is true for the other pairs. Hence it suffices to prove (3) for a set of seven points $(a-a_i) \cup x \cup y$, where $a_i$ is a
point of $a$ for some $i = 1 \dots 6$, $x = b_1$ and $y = b_2$. Let $H_i$ be the 4-space in $\mathbf P^5_a$ spanned by $a - a_i$, then (3) is satisfied if $H_i \cap < x y > = \emptyset$. Now $<x y>$ is a 
general bisecant line to $\mathbb Y_u$ because $(b_1, b_2)$ is general. So it is obvious that such a
bisecant does not intersect $H_i$. \endremark \rm \par \noindent
 Equivalently $\mathbb Y$ is the family of triples
$
(a, b_1, b_2) \in \mathbb A \times (\mathbf P^2 \times \mathbf P^2)^2
$
such that: \par \noindent
$\circ$ the linear span $<a \cup \ b_1 \ \cup b_2>$ is 6-dimensional and transversal to $\mathbf P^2 \times \mathbf P^2$, \vfill \eject  \noindent
$\circ$ the above condition (*) is satisfied. \par \noindent
\definition Let $y = (a, b_1, b_2) \in \mathbb Y$, then  $ \mathbf P^6_y := < a \ \cup \ b_1 \ \cup \ b_2>$ and $\mathbf P^{3+}_y := \mathbf P^6_y \cap \mathbf P^+$. 
\enddefinition \rm \par \noindent   
To continue with our elementary constructions we need now to consider the  Hilbert scheme of conics in $\mathbf P^2 \times \mathbf P^2$. This is split in 3 irreducible  components, according to the bidegree of a conic with respect to  first Chern classes of ${\cal O}_{\mathbf P^2 \times \mathbf P^2}(1,0)$ and ${\cal O}_{\mathbf P^2 \times \mathbf P^2}(0,1)$. We will be interested only in the component of bidegree $(1,1)$, that is, in conics $B$ such that $\pi_1(B)$ and $\pi_2(B)$ are lines. \par \noindent
\definition $\mathbb B$ is the Hilbert scheme of smooth conics of bidegree $(1,1)$ in $\mathbf P^2 \times \mathbf P^2$.
\enddefinition \rm \par \noindent  
\proposition $ \mathbb B$ is rational.
\endproposition
\proof If $B \in \mathbb B$ then $B \subset L_1 \times L_2 \subset \mathbf P^2 \times \mathbf P^2 \subset \mathbf P^8$,  where $L_i$ is the line $\pi_i(B)$ and $ i=1,2$. Note that $L_1 \times L_2$ is  embedded in $\mathbf P^2 \times \mathbf P^2$ as a smooth quadric, so that  $B \in  \mid {\cal O}_{L_1 \times L_2}(1) \mid$ $ = \mathbf P^3$.  On the other hand $(L_1,L_2)$ is a point of $\mathbf P^{2*} \times \mathbf P^{2*}$. It is easy to conclude that $\mathbb B$ is biregular to an open set of the $\mathbf P^3$-bundle on $\mathbf P^{2*} \times \mathbf P^{2*}$, whose fibre at $(L_1,L_2)$ is $ \mid {\cal O}_{L_1 \times L_2}(1) \mid$. \endproof \rm \par \noindent
\lemma Let $b_1, b_2$ be distinct  points in $\mathbf P^2 \times \mathbf P^2$ such that their projections 
$\pi_i(b_1)$ and $\pi_i(b_2)$ are distinct points of $\mathbf P^2$, for $i = 1,2$. Then the family of conics
$$
\lbrace B \in \mathbb B \ / \ b_1, b_2 \in B \rbrace
$$
is a pencil of plane sections of a smooth quadric surface in $\mathbf P^2 \times \mathbf P^2$.
\endlemma 
\proof Let $B \in \mathbb B$ be a conic containing $b_1, b_2$. Then $L_i := \pi_i(B)$ is the line joining
$\pi_i(b_1)$ to $\pi_i(b_2)$, for $i = 1,2$. Moreover it is obvious that $B \subset L_1 \times L_2$ and
that $L_1 \times L_2$ is embedded in $\mathbf P^2 \times \mathbf P^2$ as a smooth quadric surface.
This implies the statement. \endproof \rm \par \noindent
As a further step we introduce a suitable family of pairs $(y,B) \in \mathbb Y \times \mathbb B$.
\definition $\mathbb O$ is the family of pairs $(y, B) \in  \mathbb Y \times \mathbb B$ such that: \par \noindent
$\circ$ $y = (a, b_1, b_2)$ and $\lbrace b_1, \ b_2 \rbrace \subset B$, \par \noindent
$\circ$  the linear span of   $<\mathbf P^6_y \cup B>$ is a hyperplane in $\mathbf P^8$ transversal to $\mathbf P^2 \times \mathbf P^2$.
\enddefinition \rm \par \noindent
\proposition $\mathbb O$ is rational. \endproposition
\proof Let $p: \mathbb O \to \mathbb Y$ be the projection map, $ \forall \ y = (a, b_1, b_2)$ we have that $p^*(y)$
is open in $\mathbf P^1$. Indeed let $\pi_i(B) = L_i$, ($i = 1,2$),  then $p^*(y)$ is the family of pairs $(a,B)$ such that $B$ is a smooth plane section of $L_1 \times L_2$ passing through $b_1$, $b_2$. It  is standard to conclude that then $p: \mathbb O \to \mathbb Y$ realizes $\mathbb O$ as an open subset
of a $\mathbf P^1$-bundle over $\mathbb Y$. Hence $\mathbb O$ is rational. \endproof \rm \par \noindent
 \rm \par \noindent
\definition Let $o = (y,B) \in \mathbb O$, we denote  the hyperplane of $\mathbf P^8$ spanned by  
$\mathbf P^6_y \cup B$ as
$$
\mathbf P^7_o,
$$
moreover we denote  the 4-dimensional intersection of \ $\mathbf P^7_o$ \  with \ $\mathbf P^+$ \ as 
\ $\mathbf P^{4+}_o$.
\enddefinition \rm \par \noindent 
\remark \rm Let us point out that, since $\mathbf P^7_o$ contains  $ \mathbf P^-$, its equation is $\iota^*$-invariant. \endremark \rm \par \noindent The fixed set of $\iota / \mathbf P^7_o$ is $\mathbf P^- \cup \mathbf P^{4+}_o$, let us consider  in addition the threefolds
$$
\tilde T_o \  := \  \mathbf P^7_o \cap \mathbf P^2 \times \mathbf P^2 \ \text {and} \ T_o \ := \ \mathbf P^{4+}_o \cap \Sigma.
$$
$\tilde T_o$ is a smooth 3-fold whose hyperplane sections are sextic Del Pezzo surfaces.  Moreover $\iota$ acts on $\tilde T_o$ and  the set of fixed points of $\iota / \tilde T_o$ is the smooth rational normal quartic curve
$$
R_oÊ:= \Delta \cap \tilde T_o \subset \mathbf P^{4+}_o
$$
 It is clear that $T_o$ is the quotient of $\tilde T_o$ by $\iota/ \tilde T_o$, in particular  we have the commutative diagram
$$
\CD
{\mathbf P^-} @<p<< {\mathbf P^7_o} @>q>> {\mathbf P^{4+}_o}  \\
{\cup} @. {\cup} @. {\cup} \\
{\mathbf P^-} @<{p/ \tilde T_o }<< {\tilde T_o} @>{q/ \tilde T_o} >> {T_o} \\
 \endCD
$$ 
\proposition The cubic 3-fold $T_o$ is the variety of bisecant lines to $R_o \subset \mathbf P^{4+}_o$. \endproposition
\proof Since $\Delta = Sing \ \Sigma$  it follows that  $\mathbf P^{4+}_o\cap \Sigma$ is the secant variety of $\Delta \cap \mathbf P^{4+}_o$.
\endproof \rm \par \noindent

\section {K3 surfaces and marking 0-cycles in  $\mathbf P^2 \times \mathbf P^2$}
Let $o \in \mathbb O$ then a  smooth quadratic section of $\tilde T_o$ is a K3 surface. The K3 surfaces we want form a special family of such quadratic sections. To construct it let us define the universal divisor
$$
\tilde {\cal T} = \lbrace (o,x) \in \mathbb O \times \mathbf P^2 \times \mathbf P^2 \ / \ x \in \tilde T_o \rbrace
$$
and the universal cycle
$$
{\cal Z} = \lbrace (o,x) \in  \tilde {\cal T}  \ / \ x \in \ a \cup B, \ \text {where $o = (y, B)$ with $y = (a,b_1,b_2)$} \rbrace
$$
over the parameter space $ \mathbb O$. Let  $ \alpha: \tilde {\cal T} \to  \mathbb O \ \text {and} \ \beta: {\tilde {\cal T}} \to \mathbf P^2 \times \mathbf P^2 $ be the natural projections. Then, at the point $o = (a,B)$,  the fibre of $\alpha$  is $\tilde T_o$ and the fibre of $\alpha / {\cal Z}$ is
$
Z_o \  := \ a \cup B.
$
We consider the ideal sheaf ${\cal I}_{{\cal Z} / {\cal T}}$ of $\cal {Z}$ and then the sheaf
$$
\textcolor{Red}{{\cal V}} :=  \alpha_* ({\cal I}_{{\cal Z} / {\cal T}} \otimes \beta^*{\cal O}_{\mathbf P^2 \times \mathbf P^2}(2)).
$$
At the point $o \in {\mathbb O}$ the  fibre of $\textcolor{Red}{{\cal V}}$  is  the vector space
$$
H^0( {\cal I}_{Z_o/ \tilde T_o}(2)).
$$
Since $o$ is general the dimension of such a vector space is constant. Hence $\textcolor{Red}{{\cal V}}$ is a vector bundle. The involution $\iota^*$ acts on each fibre of $\textcolor{Red}{{\cal V}}$. So we have  in $ \textcolor{Red}{{\cal V}}$ a natural subbundle  whose fibre at $o$ is 
$$
H^0( {\cal I}_{Z_o/ \tilde T_o}(2))^+,
$$
i.e. the $+1$ eigenspace of $\iota^*_o$. The projectivization of such a subbundle will be denoted as
$$
\tilde {\mathbb S}.
$$
If $o \in {\mathbb O}$ then  $\tilde {\mathbb S}_o =  \mid {\cal I}_{Z_o / \tilde T_o}(2)) \mid^+$: we describe some properties of such a linear system.
\proposition If $o \in \mathbb O$ and $\tilde S \in \tilde {\mathbb S}_o$ are general then  $\tilde S$  is a smooth K3 surface. Furthermore: \par \noindent
(1) $\tilde S$ is the section of \ $\tilde T_o$ \  by $q^* Q$, for some $Q \in \mid {\cal O}_{\mathbf P^+}(2) \mid $, \par \noindent
(2) $\iota/\tilde S$ is an involution with exactly eight fixed points, which are the points of $\tilde S \cap \Delta$, \par \noindent
(3) let $S$ be the quotient of $\tilde S$ by $\iota/ \tilde S$ then $S = q(\tilde S)$,  moreover $S = \mathbf P^{4+}_o \cap \Sigma \cap Q$, \par \noindent
(4) let $y = (a, b_1, b_2)$ and let ${\cal I}_{a/\tilde S}$ be the ideal sheaf of $a$ in $\tilde S$ then the linear system
$$
\mid {\cal I}_{a / \tilde S}(1) \mid^+
$$
is a pencil and its general element is a smooth, irreducible canonical curve $\tilde F$ of genus 7. 
\endproposition
\proof  (1) \rm \ \ \  Let $\tilde T$ be a general section of $\mathbf P^2 \times \mathbf P^2$ by a hyperplane through $\mathbf P^-$, then $\tilde T$ is smooth and  contains a smooth conic $ B \in \mathbb B$. Let ${\cal J}$ be the ideal sheaf of $B$ in $\tilde T$, since $\tilde T$ is generated by quadrics it follows that the base locus of
$ \mid {\cal J}(2) \mid$ is exactly $B$ and that the base locus of $\mid {\cal J}(2) \mid^+$ is $B \cup \iota(B)$. This easily implies the following property: $\forall \ x \in B \cup \iota(B)$  the set of the elements of $\mid {\cal J}(2)\mid^+$ which are singular at $x$ has codimension 2. Then the locus of the elements which are singular at some $x \in B$ has codimension 1. Hence a general $\tilde S \in \mid {\cal J}(2) \mid^+$ is smooth along $B \cup \iota(B)$ and then it is smooth by Bertini's theorem. Such a $\tilde S$ is a K3 surface and  a quadratic section of $\tilde T$ defined by by a $\iota^*$-invariant equation. Hence $\tilde S = \tilde T_o \cap q^*Q$, for some $Q \in \mid {\cal O}_{\mathbf P^+}(2) \mid$. \ \ (2) \ \  \rm It follows from (1)
that $\iota/\tilde S$ is an involution on $\tilde S$ whose set of fixed points is 
$
\Delta \cap \tilde T \cap Q.
$
\ \ (3) \ \ Finally the surface $S := q(\tilde S)$ is the quotient of $\tilde S$ by $\iota/ \tilde S$ and it is a sextic K3 surface. More precisely $S = \mathbf P^4 \cap \Sigma \cap Q$, where $\mathbf P^4$ is the 4-dimensional linear space spanned by $q(\tilde S)$ in $\mathbf P^+$. \ \ (4) \ \ \rm Note that, for a general $\tilde T$
as above, we certainly have $\tilde T = \tilde T_o$ for some $o = (y,B)$, where $B$ is as above, $y = (a, b_1, b_2)$ and $a \in \mathbb A_i,$. This implies the statement. \endproof \rm \par \noindent
The proposition describes a  well known situation: by definition a Nikulin involution $i$ on a K3 surface  $\tilde X$ is an involution with exactly 8 fixed points. In particular $X := \tilde X/<i>$ is a K3 surface  and the image of the set of fixed points by the quotient  map is an even set of 8 nodes on $X$. 
\remark \rm Let $V$ be the family of all K3 surfaces $\tilde X \subset \mathbf P^2 \times \mathbf P^2$
which are complete intersections of a hyperplane and a quadratic section and such that  $\iota(\tilde X) = \tilde X$. $PGL(3)$ admits a diagonal action on $\mathbf P^2 \times \mathbf P^2$ and hence on $V$: it turns out that its GIT-quotient is the moduli space of K3 surfaces endowed with a Nikulin involution and with the genus two polarization $(\pi_1/\tilde X)^*{\cal O}_{\mathbf P^2}(1)$.
\endremark \rm \par \noindent
Fix $o = (y,B) \in \mathbb O$ and $\tilde S \in \tilde {\mathbb S}_o$ then $\mathbf P^6_y$ is a hyperplane in the ambient space $\mathbf P^7_o$ of $\tilde S$. Let
$$
\tilde F:= \mathbf P^6_y \cap \tilde S
$$
then the equation of $\tilde F$ is $\iota^*$-invariant and $\iota / \tilde F$ is a fixed-point-free involution on $\tilde F$. The curve 
$$
F := q(\tilde F) \subset \Sigma \cap \mathbf P_y^{3+}
$$
is the quotient of $\tilde F$ by $\iota / \tilde F$ and the quotient map is the \'etale double covering
$$
q / \tilde F: \tilde F \to F.
$$
\proposition For general $o = (y,B) \in \mathbb O$ and $\tilde S \in \tilde {\mathbb S}_o$ we have: \par \noindent
$\circ$ $\tilde F$ is a smooth canonical curve of genus $7$ in $\mathbf P^6_y$,  \par \noindent
$\circ$ $\tilde F$ is endowed with the fixed point free involution $\iota / \tilde F$, \par \noindent
$\circ$ $F$ is a smooth, canonical curve of genus 4 in $\mathbf P^{3+}_y$.
\endproposition
\proof Let $\tilde F$ be a general section of $\tilde S$ by a hyperplane containing $\mathbf P^-$. Then $\tilde F$ is a smooth canonical curve of genus 7, $\iota / \tilde F$ is a fixed-point-free involution and $
q/\tilde F: \tilde F \to F = q(\tilde F)$ is its quotient map.  As in section 2 we consider  the Norm map $Nm: Pic^6(\tilde F) \to Pic^6(F)$ and the two connected components $P^+$ and $P^-$ of $Nm^{-1}(\omega_F)$. For a general $N \in P^-$ we have $N = {\cal O}_{\tilde F}(a')$, where $a'$ is a smooth
divisor, $q: a' \to q(a')$ is bijective and $h^0(N) = 1$. The latter equality implies that the 6 points of
$a'$ are linearly independent.  Let $z \in a'$ and let $a'' = a' - z + \iota(z)$, then $h^0({\cal O}_{\tilde F}(a'')) = 2$. Indeed it follows from prop. 2.1 that $N(i(z)-z)) \in P^+$, that is, $h^0(N(i(z)-z))$ is even. 
On the other hand $a''$ is effective and $h^0({\cal O}_{\tilde F}(a'')) \leq 2$, because $a''-i(z)$ consists of 5 linearly independent points. Hence $h^0({\cal O}_{\tilde F}(a'')) = 2$. \\
 Applying geometric Riemann-Roch to the canonical curve $\tilde F$, the linear span $\Lambda = < a''>$ is a 4-space, moreover $\Lambda$ intersects $\mathbf P^-$ along a line. It is easy to conclude that the pair $(a'',z)$, is a point of  $\mathbb A$.
Finally let $y' = (a', b_1, b_2)$ where $b_1 + b_2 = B \cdot \tilde F$, then $o' = (y', B)$ is a point of
$\mathbb O$ and $\tilde S \in \tilde {\mathbb S}_{o'}$. The statement clearly holds true  for $o'$. Hence it holds on an open dense subset of $\tilde {\mathbb S}$.
\endproof \rm \par \noindent
Keeping the above notation we have the commutative diagram
 $$
\CD
{\mathbf P^-} @<p<< {\mathbf P^7_o} @>q>> {\mathbf P^{4+}_o}  \\
{\cup} @. {\cup} @. {\cup} \\
{\mathbf P^-} @<{p/ \tilde S}<< {\tilde S} @>{q/ \tilde S} >> {S} \\
{\cup} @. {\cup} @. {\cup} \\
{\mathbf P^-} @<{p/ \tilde F}<< {\tilde F} @>{q/ \tilde F} >> {F} \\
\endCD
$$  \\
We are now in position to define the main family of $0$-cycles of this paper: in particular these are subschemes $d \subset \mathbf P^2 \times \mathbf P^2$ of length 8 such that $q_*d$ is the base locus of a net of quadrics of $\mathbf P^{3+}_y$. So, for a general $d$, $q_*d$ is a hyperplane section of a smooth canonical curve in $\mathbf P^{4+}_o$. 
Consider a general $o = (y, B) \in \mathbb O$ and a general $\tilde S \in \tilde {\mathbb S}_o$, then
consider the curves $\tilde F = \mathbf P^6_y \cap \tilde S$ and $F = q(\tilde F)$. Let  $y = (a, b_1, b_2)$
then one has on $\tilde F$  the smooth divisors of degree 8
$$
m_1 := a + b_1 + \iota (b_2) \ \text { and } \ m_2 = a + b_1 + b_2.
$$ 
\definition A marking $0$-cycle of type $i = 1$ or $2$ is a triple $(o, \tilde S, d)$ such that  
$$
d \in \mid {\cal O}_{\tilde F}(m_i) \mid.
$$
The family of all marking 0-cycles of type $i$ will be denoted as
$$
\mathbb D_i.
$$
\enddefinition \rm \vfill \eject  \noindent
\proposition $\mid {\cal O}_{\tilde F}(m_i) \mid$ is a base-point-free pencil for $i = 1,2$.
\endproposition
\proof By the definition of $\mathbb Y$, $m_i$ is not in a hyperplane. Since $\tilde F$
is embedded in $\mathbf P^6_y$ as a canonical curve,  geometric Riemann-Roch implies 
$h^1({\cal O}_{\tilde F}(m_i)) = 0$. Hence $h^0({\cal O}_{\tilde F}(m_i)) = 2$
and  $\mid m_i \mid$  is a pencil.  This is base-point-free if the points of  any divisor of degree 7 contained in $m_i$ are linearly independent. Since $y = (a, b_1, b_2) \in \mathbb Y$, the latter  condition is satisfied by the definition of $\mathbb Y$.
\endproof \rm \par \noindent
Consider the projection map $ {\mathbb D}_i \to \tilde {\mathbb S}$ and observe that its fibre at
$(o, \tilde S)$ is $\mid {\cal O}_{\tilde F}(m_i) \mid = \mathbf P^1$. This, with some more standard work we omit for brevity, implies that
\proposition ${\mathbb D}_1$ and $\mathbb D_2$ are  $\mathbf P^1$-bundles over $\tilde {\mathbb S}$, hence they are rational.
\endproposition \rm \par \noindent
\section {The family of marked curves $\tilde C$}
Since we have all the ingredients, it is time to cook them up to construct a family of pairs
$$
(\pi, d)
$$
where $\pi: \tilde C \to C$ is an \'etale double covering and $\tilde C$ is a curve of genus 9, marked by a divisor  $d$ satisfying the following condition:  \  $\pi_*d \in \mid \omega_C \mid$ and $h^0({\cal O}_{\tilde C}(d)) = 1$ or $h^0({\cal O}_{\tilde C}(d)) = 2$.  \par \noindent Let $P(\pi)$ be the Prym of $\pi$:  in the
former case ${\cal O}_{\tilde C}(d)$ is a point of the model $P^-$ of $P(\pi)$,  in the latter 
${\cal O}_{\tilde C}(d)$ is a point of the theta divisor $\Xi$ of $P(\pi)$, embedded in its model $P^+$
(see section 2). \par \noindent
Once more we will consider a general point $o = (y,B) \in \mathbb O$ and a general surface $\tilde S$
in the linear system $\tilde {\mathbb S}_o$. As usual let $\tilde F = \mathbf P^6_y \cap \tilde S$  then
we have on $\tilde S$ the linear system
$$
\mid \tilde F + B + \iota^*B \mid.
$$
\proposition $\mid \tilde F + B + i^*B \mid$ is  a 9-dimensional and base-point-free linear system of smooth irreducible curves of genus 9. \endproposition
\proof Let $H :=  \tilde F + B + i^*B $. Since $\tilde F$ is base-point-free, $\mid H \mid$ is base-point-free on the open set $\tilde S - (B \cup i^*B)$. Note that ${\cal O}_{B \cup i^*B}(H)$  is the trivial sheaf ${\cal O}_{ B \cup i^*B}$. Then $\mid H \mid$ is base-point-free if the
restriction map $H^0({\cal O}_{\tilde S}(H)) \to H^0({\cal O}_{B \cup i^*B}(H))$ is surjective. This follows
from the vanishing of $H^1({\cal O}_{\tilde S}(\tilde F))$ and the long exact sequence associated to the standard exact sequence
$$
0 \to {\cal O}_{\tilde S}(\tilde F) \to {\cal O}_{\tilde S}(H) \to {\cal O}_{B \cup i^*B}(H) \to 0.
$$
Since $H^2 > 0$ and $\mid H \mid$ is base-point-free, a general $\tilde C \in \mid H \mid$ is a smooth, irreducible curve of genus 9, moreover $h^1({\cal O}_{\tilde S}(H)) = 0$, (see [SD]). In particular it  follows  $dim \mid H \mid \\  = 9$.  \endproof 
\rm \par \noindent
Observe that $\mid \tilde F + B + \iota^*B \mid$ is invariant under the action of the involution $\iota$. In this linear system we have therefore the projectivized eigenspace $\mid \tilde F + B + \iota^*B \mid^+$ of $\iota$.
\proposition A general $\tilde C \in \mid F + B + i^*B \mid^+$ is a smooth, irreducible curve of genus 9,
endowed with a fixed-point-free involution $i: \tilde C \to \tilde C$.
\endproposition \vfill \eject \noindent
\proof  With the usual notations let $F = q(\tilde F)$ and $S = q(\tilde S)$,  then $F$ is the section of $S  \subset \mathbf P^{4+}_y$ by the hyperplane $\mathbf P^{3+}_y$. We know that
$
S = \mathbf P^{4+}_y \cap \Sigma \cap Q
$
where $Q \subset \mathbf P^{4+}_y$ is a quadric hypersurface. In addition we also know that $S$ is a K3 surface, singular exactly at the 8 nodes of $\Delta \cap \Sigma$. 
Since $ \mid \tilde F + B + i^*B \mid^+$ is the projectivized $+1$ eigenspace of $\iota / \tilde S$, it follows that
$
\mid \tilde F + B + i^* B \mid^+ = (q/\tilde S)^* \mid F + \overline B \mid,
$
where $\overline B$ is the conic $q_*B$. Since each element of $\mid \tilde F + B + \iota^*B \mid$ is connected, the proposition follows if we show that a general $C \in \mid F + \overline B \mid$ is smooth, irreducible. To prove such a property just apply the same proof used in the previous proposition.  \endproof \rm \par \noindent
\proposition Let $\tilde C$  be a smooth element of  $\mid \tilde F + B + i^* B \mid^+$ and let
${\cal I}_{\tilde C / \tilde T}$ be the ideal sheaf of $\tilde C$ in the threefold $\tilde T = < \tilde S > \cap \ \mathbf P^2 \times \mathbf P^2$, then: \par \noindent
(1) ${\cal I}_{\tilde C / \tilde T}(2)$ is acyclic and $h^0({\cal I}_{\tilde C / \tilde T}(2)) = 3$, \par \noindent
(2) $h^0({\cal I}_{\tilde C / \tilde T}(2))^- = 0$.
\endproposition
\proof (1) Consider the standard exact sequence of ideal sheaves
$$
0 \to {\cal I}_{\tilde C / \tilde S}(2) \to {\cal I}_{\tilde C / \tilde T}(2) \to {\cal I}_{\tilde S / \tilde T}(2) \to 0.
$$
We have ${\cal I}_{\tilde C / \tilde S}(2) = {\cal O}_{\tilde S}(E)$ where  $E := F - B - \iota^*B$. Note
that $h^0({\cal O}_{\tilde S}(E)) = 2$. This follows because we assume $B$ is  general in its family of conics so that $< B > \cap < \iota^* B > $ $=$ $ \emptyset$. Then $< B \cup \iota^*B >$ is a 5-space  in  $\mathbf P^7 = < \tilde S >$ and this implies $h^0({\cal O}_{\tilde S}(E)) = 2$.  Since $E^ 2 = 0$, Riemann-Roch implies that ${\cal O}_{\tilde S}(E)$ is acyclic. On the other hand ${\cal I}_{\tilde S / \tilde T}(2)$
is just  ${\cal O}_{\tilde T}$, which is acyclic. Then ${\cal I}_{\tilde C / \tilde T}$ is acyclic too and  the associated long exact sequence yelds $h^0({\cal I}_{\tilde C}(2)) = 3$. \\
(2) Note that $C = q(\tilde C)$ is a canonical curve of genus 5 in the space $\mathbf P^{4+}_y$ previously considered. In particular observe that $T = q(\tilde T)$ is a cubic hypersurface in $\mathbf P^{4+}_y$ and $h^0({\cal I}_{C / T}(2)) =3$, where ${\cal I}_{C / T}$ is the ideal sheaf of $C$ in $T$. Then (1) implies that $H^0({\cal I}_{\tilde C /  \tilde T}(2)) = q^*H^0({\cal I}_{C / T}(2))$. In particular
each section $s \in H^0({	\cal I}_{\tilde C / \tilde T}(2))$ is $\iota^*$-invariant, that is $h^0({\cal I}_{\tilde C / \tilde T}(2))^- = 0$. \endproof \rm \par \noindent
 \lemma[\it Parity lemma] Let $d \subset \tilde S$ be a smooth scheme of length 8 such that: \\
(1) $ \tilde F = < d > \cap \  \tilde S$ is a smooth hyperplane section transversal to $B + \iota^*B$. \\
(2)  $d \sim  a + b' + b''$ where $\lbrace b', b'' \rbrace \subset B \cup \iota^*B$, $a$ is effective and $a \cap (B \cup \iota^*B) = \emptyset $. \\
(3) $\mid a + b' + b'' \mid$ is a base-point-free pencil and $d$ is general in it. \\
(4) A general element of $ \mid a + b' + b'' \mid$ is included in a smooth $\tilde C \in  \mid \tilde F + B + \iota^* B \mid$. \\
Then, for the ideal sheaf of $d$ in $\tilde S$, we have: \par \noindent
$\circ$ $h^0({\cal I}_{d / \tilde S}(\tilde C)) = 2$ \  if $\# ((b' + b'') \cap B)$ is odd ,   \par \noindent
$\circ$ $ h^0({\cal I}_{d / \tilde S}(\tilde C)) = 3$ \ if $\# ((b' + b'') \cap B)$ is even.
 \endlemma
\proof   If $\# ((b' + b'') \cap B)$ is  odd we consider the standard exact sequence of sheaves
$$
0 \to {\cal I}_{d / \tilde S}(\tilde F) \to {\cal I}_{ d / \tilde S}(\tilde C) \to {\cal J}_d \to 0,
$$
where the map ${\cal I}_{d / \tilde S}(\tilde F) \to {\cal I}_{ d / \tilde S}(\tilde C)$ is the natural inclusion of ideal sheaves and ${\cal J}_d$ its cokernel. By (1) $d$ spans a  hyperplane in the 7-space $< \tilde S >$, hence $h^0({\cal I}_{d / \tilde S}(\tilde F) ) = 1$ . Since $deg \ d = 8$ we have also $h^1({\cal I}_{d / \tilde S} (\tilde F)) = 1$. Since  $\tilde C B = \tilde C \iota^*B = 0$ we have ${\cal J}_d = {\cal O}_{B \cup \iota^*B} $ if $d \cap B \cup \iota^*B = \emptyset$. By (3) this happens for a general $d$. Passing to the associated long exact sequence it follows 
$$
h^0({\cal I}_{ d / \tilde S}(\tilde C) \geq 2, \ \ \ d \in \mid a + b' + b'' \mid.
$$
Then, by semicontinuity, the equality follows for a general $d$ if it holds for $a + b' + b'' $. To prove that
$ h^0({\cal I}_{ a + b' + b'' / \tilde S}(\tilde C) = 2 $ we remark that $B + \iota^*B$ is a fixed component of
$ \mid  {\cal I}_{a + b' + b'' / \tilde S}(\tilde C) \mid$. 
This is clear because $\tilde C$ intersects both $B$ and $\iota^*B$, while we have $\tilde C B = \tilde C \iota^*B = 0$. Therefore it follows  $h^0({\cal I}_{a + b'  + b'' / \tilde S}(\tilde C))$ 
$=$ $h^0({\cal I}_{a / \tilde S}(\tilde F))$. Now $h^0{\cal O}_{\tilde F}(a)) = 1$  because the pencil
$\mid a + b' + b'' \mid$ is base-point-free, hence $h^0({\cal I}_{a / \tilde F}(\tilde F)) = 1$ that is $h^0({\cal I}_{a / \tilde S}(\tilde F)) = 2$. \\
 If $\# ((b' + b'') \cap B)$ is  even we can assume $b' + b''  \subset B$. We consider the standard exact  sequence
$$
0 \to {\cal I}_{d / \tilde S}(\tilde F + \iota^*B) \to {\cal I}_{ d / \tilde S}(\tilde C) \to {\cal J}_d \to 0
$$
and claim that $h^0({\cal I}_{d / \tilde S}(\tilde F + B)) = 2$. This implies that $\iota^*B$ is a fixed component of the pencil $\iota^*B+\mid {\cal I}_{d / \tilde S}(\tilde F + B) \mid$ which is contained in  $ \mid {\cal I}_{ d / \tilde S}(\tilde C) \mid$. Then assumption (4) implies
$$
h^0({\cal I}_{ d / \tilde S}(\tilde C)) \geq 3, \ \ \ d \in \mid a + b' + b'' \mid.
$$
By semicontinuity the equality follows for a general $d$ if it holds for $d = a + b' + b''$. Since $b' + b''$
is in $B$ it follows, arguing as in the previous part of the proof, that $B$ is a fixed component of
$\mid {\cal I}_{d / \tilde S}(\tilde C) \mid$ and that  $h^0({\cal I}_{d / \tilde S}(\tilde C))$ $=$ $h^0({\cal I}_{a / \tilde S}(\tilde F + \iota^*B))$. Under assumption (2) we have also $a \cap \iota^*B = \emptyset$. Therefore we have  the following standard exact sequence of ideal sheaves
$$
0 \to {\cal I}_{a / \tilde S}(\tilde F) \to {\cal I}_{a / \tilde S}(F + \iota^*B) \to {\cal O}_{\iota^*B} \to 0.
$$
As above $h^0({\cal O}_{\tilde F}(a)) = 1$ so that $h^0({\cal I}_{a / \tilde S}(\tilde F)) =2$ and $h^1({\cal I}_{a / \tilde S}(\tilde F)) = 0$. Then it follows that $h^0({\cal I}_{a / \tilde S}(\tilde F + \iota^*B)) = 3$. \\
To complete the proof we prove our previous claim that $h^0({\cal I}_{d / \tilde S}(\tilde F + B)) = 2$: \ \ 
since $\tilde F$ is very ample, the map defined by $\mid \tilde F + B \mid$ is a birational morphism 
$\sigma: \tilde F \to \tilde F_o$ onto its image $\tilde F_o$ and $o = \sigma(B)$ is a node, (cfr. [SD], see also [S] lemma 2.4).
Since $\sigma(b') = \sigma(b'') = o$, ${\cal O}_{\tilde F}(a + b' + b'')$ descends to a line bundle $E$ on $\tilde F_o$ such that $h^0(E) = h^0({\cal O}_{\tilde F}(a + b' + b'')) = 2$. By
Riemann-Roch on $\tilde F_o$, it follows $h^0(\omega_{\tilde F_o} \otimes E^{-1}) = 1$.  This easily implies that $h^0({\cal I}_{a + b' + b'' / \tilde F}(\tilde F + B)) = 2$. \endproof
 \rm \par \noindent
Keeping the previous notations let us consider now the linear system defined as follows:
$$
\mid {\cal I}_{d / \tilde S}(\tilde C) \mid^+ \ := \mid {\cal I}_{d / \tilde S}(\tilde C) \mid \cap \  \iota^* \mid {\cal I}_{d / \tilde S}(\tilde C) \mid.
$$
This consists of  curves $\tilde C$ containing $d$ and having a $\iota^*$-invariant equation. We want to apply the parity lemma  to marking $0$-cycles:
\proposition Let $(o, \tilde S,d)$ be a general marking $0$-cycle of type $i$ $=$ 1 or 2, then we have: \\
(1) $\mid {\cal I}_{d / \tilde S}(\tilde C) \mid^+$ is a pencil. \\
(2) A general $\tilde C \in \mid {\cal I}_{d / \tilde S}(\tilde C) \mid^+$ is a smooth, irreducible curve. \\
(3) The dimension of $\mid {\cal I}_{d / \tilde S}(\tilde C) \mid$ is $i$.
 \\ 
\endproposition
\proof   Let  $\overline B = q(B)$ and $F = q(\tilde F)$ then we have
$$
 \mid ({\cal I}_{d/\tilde S}(\tilde F + B + i^*B)\mid ^+ \subset \mid \tilde F + B + i^* B \mid^+ := (q/\tilde S)^* \mid F + \overline B \mid.
$$
Since $(o, \tilde S, d)$ is a marking 0-cycle we have: $d \sim a + b' + b''$, where $a \in \mathbb A$ and $q_*a \in \mid \omega_F \mid$. Moreover $b' \neq \iota(b'')$ so that $q_*(b'+b'') = F \cdot \overline B$.
This implies that
$$
q_*d \in \mid {\cal O}_F(F + \overline B) \mid.
$$
Consider the standard exact sequence of ideal sheaves
$$
0 \to {\cal I}_{F / S}(F + \overline B) \to {\cal I}_{q_*d / S}(F + \overline B) \to {\cal I}_{q_*d / F}(F + \overline B) \to 0.
$$
By construction $B \cap Sing \ S = \emptyset$, so that ${\cal I}_{F / S}(F + \overline B)$ is the line bundle ${\cal O}_{S}( \overline B)$. Since $ \overline B^2 = -2$ and $B$ is effective, Riemann-Roch and Serre duality imply $h^1({\cal O}_{\tilde S}( \overline B)) = 0$. On the other hand $q_*d \in \mid \omega_F \mid$  implies that $ {\cal I}_{q_*d / F} (F + \overline B) \cong {\cal O}_F$. Passing to the associated long exact sequence it follows $h^0({\cal I}_{q_*d / S}(F + \overline B)) = 2$. Hence $P := \mid {\cal I}_{q_*d / S}(F + \overline B)   \mid$ is a pencil as well as
$$
q^*P = \mid {\cal I}_{d / \tilde S}(\tilde C) \mid^+.
$$
Let us show that a general $C \in P$ is smooth: since $F + \overline B \in P$ it follows that $P$ has no
fixed components. Indeed $F$ is not a fixed component because $B$ is isolated. If $B$ is a fixed
component then the moving part of $P$ is an irreducible pencil of hyperplane sections of $S$, with base locus a scheme of length 6. This implies that two points of  $q_*d$ are in $\overline B$, which is not true for a general $d$. Since $C^2 = 8$ we conclude that $q_*d$ is the base locus of $P$. Since $q_*d$ is smooth for a general $d$, a general $C \in P$ is smooth. This implies (1) and (2). (3) follows from parity lemma. \endproof \rm \par \noindent
For a general $ \tilde C \in \mid {\cal I}_{d / \tilde S}(\tilde F + B + \iota^*B) \mid$ as above $\iota / \tilde C$ is a fixed-point-free involution. So  $\tilde C$ is a smooth curve of genus 9 endowed with $\iota /  \tilde C$ and marked by $d$. Let $C = \tilde C / < \iota >$ and let
$$
\pi: \tilde C \to C
$$
be the quotient map. Then $\pi = q/ \tilde C$ and $C = q(\tilde C)$ is a canonical curve of genus 5 in  the 4-space $\mathbf P^{4+}_y$, (here $o = (y, B)$ and $(o, \tilde S,d)$ is the marking $0$-cycle considered in the previous statement).
\proposition Let $d$ and $\pi: \tilde C \to C$ be as above, then $\pi_*d \in \mid \omega_C \mid$ and
$h^0({\cal O}_{\tilde C}(d)) = i$, where $i = 1,2$ is the type of the marking $0$-cycle $(o, \tilde S, d)$.\endproposition
\proof  To see that $\pi_*d \in \mid \omega_C \mid$  just note that $\pi_*d$ is a hyperplane section of
the canonical curve $C$. Indeed we have $\pi_*d = q_*d = q(\tilde F) \cdot C$ and $q(\tilde F)$ is a
hyperplane section of $S = q(\tilde S)$. To see that $h^0{\cal O}_{\tilde C}(d)) = i$ observe that, by (3) of the previous proposition,   $h^0({\cal I}_{d / \tilde S}(\tilde F + B + \iota^*B)) = i + 1$. But this is exactly equivalent to $h^0({\cal O}_{\tilde C}(d)) = i$.
 \endproof \rm \par \noindent
\remark \rm  Let $P(\pi)$ be the Prym variety of $\pi: \tilde C \to C$ and let
$
Nm: Pic^8(\tilde C) \to Pic(C)
$
be the Norm map defined in section 2. Keeping the notations fixed there, the proposition implies that   
$
{\cal O}_{\tilde C}(d) \in Nm^{-1}(\omega_C) = P^+ \cup P^-.
$
Since $h^0({\cal O}_{\tilde C}(d)) = i$ with $i$ $=$  $1$ or $2$ we have in addition that
$$
{\cal O}_{\tilde C}(d) \in P^- \cup \Xi.
$$
\vfill \eject \noindent
$P^-$ and $P^+$ are copies of $P(\pi)$ and $\Xi \subset P^+$ is a copy of the theta divisor of $P(\pi)$.
So the proposition says that ${\cal O}_{\tilde C}(d)$ defines  a point of $P(\pi)$ if $i = 1$ and of its theta divisor if $i = 2$. \endremark \rm \par \noindent
\definition Let $i$ $=$ $1$ or $2$ then $\tilde {\mathbb C}_i$ is  the family of 4-tuples $(o, \tilde S, d, \tilde C)$, where $(o, \tilde S, d)$ is a marking 0-cycle of type $i$ and
$$
\tilde C \in \mid {\cal I}_{d/\tilde S}(\tilde F + B_{a,o} + \iota^* B_{a,o})\mid^+.
$$
\enddefinition \rm \par \noindent
\proposition $\tilde {\mathbb C}_1$ and $\tilde {\mathbb C}_2$ are rational varieties. \endproposition 
\proof Let $i =$ $1$ or $2$. It is clear that the natural projection map $\tilde {\mathbb C}_i \to \mathbb D_i$ is a  $\mathbf P^i$-bundle structure over $\mathbb D_i$, so we omit further details. Since $\mathbb D_i$ is
rational, it follows that $\tilde {\mathbb C}_i$ is rational. \endproof \rm \par \noindent
\section {The unirationality results}
Finally we use the constructions given in the previous sections  to prove our main theorem which says that both  the universal p.p.a.v. and the universal theta divisor over \textcolor{Red}{${\cal A}_4$} are unirational varieties. Let
$$
\tilde {\mathbb C} = \tilde {\mathbb C}_1 \cup \tilde {\mathbb C_2},
$$
moreover let 
$$
\ \textcolor{Red}{\tilde {\cal C}} \subset \tilde {\mathbb C} \times \mathbf P^2 \times \mathbf P^2  \ \rm and \  \textcolor{Red}{{\cal C}} \subset \tilde {\mathbb C}  \times \Sigma
$$
be respectively  the universal curve over $\tilde {\mathbb C}$ and its image via the product map $$ id_{\tilde {\mathbb C}} \times q: \tilde {\mathbb C} \times \mathbf P^2 \times \mathbf P^2 \to \tilde {\mathbb C} \times \Sigma. $$ Restricting this map  to $ \textcolor{Red}{\tilde {\cal C}}$ we obtain a finite double covering of $\tilde {\mathbb C}$-schemes
$$
\emph q:  \textcolor{Red}{\tilde {\cal C}} \to \textcolor{Red}{{\cal C}}.
$$
At $x = (o, \tilde S, d, \tilde C) \in \tilde {\mathbb C}$ the fibre map $\emph q_x:  \textcolor{Red}{\tilde {\cal C}}_x \to \textcolor{Red}{{\cal C}_{x}}$ is  the  \'etale double covering $\pi: \tilde C \to C$, where $C = q(\tilde C)$ and $\pi = q / \tilde C$.  Let us define  the following  line bundles associated to $x$:
$$
M_x := {\cal O}_{\tilde C}(d) \ \text {and} \  L_x := {\cal O}_{\tilde C} \otimes {\cal O}_{\mathbf P^2 \times \mathbf P}(1,0).
$$
$M_x$ has been studied in proposition 5.6.  Let $i$ be the involution defined by  $\pi$, then $i$ is exactly $\iota/\tilde C$. This implies that $i^*L_x = {\cal O}_{\tilde C} \otimes  {\cal O}_{\mathbf P^2 \times \mathbf P^2}(0,1)$ and that $L_x \otimes i^*L_x \cong {\cal O}_{\tilde C}(1)$. Since ${\cal O}_{\tilde C}(1) \cong \omega_{\tilde C}$ it follows that 
$$
L_x \in Nm^{-1}(\omega_C) = P^+ \cup P^-.
$$
Here $Nm: Pic^8(\tilde C) \to Pic^8(C)$ is the Norm map defined  in section 2. Recall also that, by definition, an element $M$ of $Nm^{-1}(\omega_C)$ is in $P^+$, (in $P^-$), iff $h^0(M)$ is even, (odd). We want to study the Prym-Petri map of $L_x$, that is the multiplication map
$$
\mu_{L_x}: [H^0(L_x) \otimes H^0(\iota^*L_x)]^- \to H^0(\omega_{\tilde C})^-.
$$
\proposition $L_x$ is an element of $P^-$, $h^0(L_x) = 3$ and the Prym-Petri map of $L_x$ is injective.
 \endproposition
\proof  First let us show that $h^0(L_x) = 3$. We put $L := L_x$, then we recall that 
$
\tilde C \ \in \ \mid \tilde F + B + \iota^*B \mid \ \text {and that} \ \tilde F \sim F_1 + F_2
$
where $ \mid \tilde F_1 \in {\cal O}_{\tilde S}(1,0) \mid $ and $ \mid \tilde F_2 \in {\cal O}_{\tilde S}(0,1) \mid $. Hence one has the standard exact sequence
$$
0 \to {\cal O}_{\tilde S}(-\tilde F_2 - B - \iota^*B) \to {\cal O}_{\tilde S}(\tilde F_1) \to L \to 0.
$$
Since $h^0({\cal O}_{\tilde S}(F_1)) = 3$, it follows $h^0(L) = 3$ if  $h^1({\cal O}_{\tilde S}(-\tilde F_2 - B - \iota^*B)) = 0$. To prove this observe that  $(F_2 + B + \iota^*B) B$ $=$
$(F_2 + B + \iota^*B)\iota^*B = -1$. Then the standard exact sequence
$$
0 \to {\cal O}_{\tilde S}(F_2) \to {\cal O}_{\tilde S}(F_2 + B + \iota^*B) \to {\cal O}_{B + \iota^*B}(F_2 + B + \iota^*B) \to 0
$$
implies $h^0({\cal O}_{\tilde S}(F_2+B+\iota^*B))$ $=$ $h^0({\cal O}_{\tilde S}(F_2)) = 3$. Applying  Riemann-Roch and Serre duality  we conclude that $h^1({\cal O}_{\tilde S}(F_2+B+\iota^*B))$ $=$  $h^1({\cal O}_{\tilde S}(-F_2-B-\iota^*B)) = 0$. \\
We are left to show that the Prym-Petri map of $L$  is injective. Consider  the Petri map
$$
\mu: H^0(\tilde L) \otimes H^0(i^* \tilde L) \to H^0(\omega_{\tilde C})
$$
then  $Im (\mu)$ defines the embedding $\tilde C \subset \mathbf P^2 \times \mathbf P^2 \subset \mathbf P^8$. This implies that  $Ker \ \mu \cong H^0({\cal I}_{\tilde C}(1))$, where ${\cal I}_{\tilde C}$ is the ideal sheaf of $\tilde C$ in $\mathbf P^8$.  Now $Ker \ \mu$ is 1-dimensional and it is generated by a $+1$ eigenvector of $i^*$. Indeed we know that $\tilde C \subset \tilde S \subset \mathbf P^7_o$, the equation of  $\mathbf P^7_o$ being $\iota^*$-invariant. Hence  there exists a non-zero $h \in H^0(L) \otimes H^0(i^*L)$ such that $i^*h = h$ and $\mu(h) = 0$. Since $h$ is a $+1$ eigenvector of $i^*$, the injectivity of the Prym-Petri map of $L$ follows if $dim \ Ker \ \mu = 1$, that is
$h^0({\cal I}_{\tilde C}(1)) = 1$. To prove this  consider the standard  exact sequence of ideal sheaves on $\mathbf P^8$
$$
0 \to {\cal I}_{\tilde S}(1) \to {\cal I}_{\tilde C}(1) \to {\cal I}_{\tilde C / \tilde S}(1) \to 0.
$$
From $h^0({\cal I}_{\tilde C / \tilde S}(1)) = h^0({\cal O}_{\tilde S}(-B - \iota^*B)) = 0$ and $h^0({\cal I}_{\tilde S}(1)) = 1$, it follows $h^0{\cal I}_{\tilde C}(1)) = 1$. 
\endproof \rm \par \noindent
$L_x$ and $M_x$ admit natural extensions to $ \textcolor{Red}{\tilde {\cal C}}$, indeed we have on $\textcolor{Red}{\tilde {\cal C}}$ the following line bundles:
$$
\textcolor{Red}{{\cal L}}Ê:= {\cal O}_{ \textcolor{Red}{\tilde {\cal C}}} \otimes \beta^* {\cal O}_{\mathbf P^2 \times \mathbf P^2}(0,1)  \ \text {and} \  {\cal M } \ := \ {\cal O}_{ \textcolor{Red}{\tilde{\cal C}}}(D),
$$
where $D = \lbrace (o, \tilde S, d, \tilde C, z)  \in \tilde {\mathbb C} \times \mathbf P^2 \times \mathbf P^2 \ / \ z \in d \rbrace$ is the universal marking $0$-cycle and  $\beta$ is the second projection of $\tilde {\mathbb C} \times (\mathbf P^2 \times \mathbf P^2)$.  For any $x \in \tilde {\mathbb C}$ the fibre $\textcolor{Red}{\tilde {\cal C}}_x$ is $\tilde C$, moreover  the restrictions of $\textcolor{Red}{{\cal L}}$ and $\textcolor{Red}{{\cal M}}$ to ${\textcolor{Red}{\tilde {\cal C}}}_x$ are respectively $L_x$ and $M_x$. Let us also define the following moduli spaces: \par \noindent
\definition Let $i \in \lbrace 1, 2 \rbrace$ then $\textcolor{Red}{{\cal W}_i}$ is the moduli space of triples $(\pi, M, L)$ such that: \par \noindent
(1) $\pi: \tilde C \to C$ is a general  \'etale double cover of a smooth, irreducible curve of genus 5, 
\par \noindent
 (2) $M$ is  a point  in $U_i := \lbrace M \in P^- \cup P^+  \ / \ h^0({\cal O}_{\tilde C}(d)) = i \ \rbrace $. 
\par \noindent
(3)  $L$ is general in $W^2(\pi) := \lbrace M \in P^- \ / \ h^0(M) \geq 3 \rbrace$ i.e. its  Prym-Petri map is injective.
\enddefinition \rm \par \noindent
\remark \rm $W^2(\pi)$ is a Prym-Brill-Noether locus: it is a scheme defined as in [W]. Counting dimensions it follows that the Prym-Petri map of $L$ cannot be injective if $h^0(L) > 3$. Therefore we are implicitly assuming $h^0(L) = 3$ in the above definition.
\endremark \rm
\remark \rm $U_1$ is open in the  the model $P^-$ of the Prym of $\pi$, more precisely $U_1 = P^- - W^2(\pi)$. A point $M \in U_2$ is just a point $M \in P^+$ such that $h^0(M) = 2$. Therefore $U_2$ is open
in the model $\Xi \subset P^+$ of the theta divisor of the Prym of $\pi$. Notice also that the restriction to
$P^-$ of the natural theta line bundle of $Pic^8(\tilde C)$ is twice  the principal polarization of the Prym, ([M2]). \endremark

\definition  $\textcolor{Red}{\textcolor{Red}{{\cal U}}_i}$ is the moduli space of $(\pi, M)$ where $\pi$ and $M \in U_i$ are as above, ($i = 1, 2$). 
\enddefinition \rm \par \noindent
For $i = 1,2$ let us consider the forfgetful map
$$p_i: \textcolor{Red}{\textcolor{Red}{{\cal U}}_i} \to \textcolor{Red}{{\cal A}_5}$$ 
onto the Prym moduli space \textcolor{Red}{$\textcolor{Red}{{\cal A}_5}$}. At the moduli point of $\pi$ the fibre $U_1$  of $p_1$  is a copy of  the  Prym of $\pi$ and the fibre $U_2$ of $p_2$ is a copy of its  theta divisor. 
\proposition (1) $\textcolor{Red}{{\cal U}}_1$ dominates the universal p.p.a.v. over \textcolor{Red}{${\cal A}_4$}.  \\
(2) $\textcolor{Red}{{\cal U}}_2$ dominates the universal theta divisor over \textcolor{Red}{${\cal A}_4$}. \endproposition
\proof Recall that the Prym map $p_5: \textcolor{Red}{{\cal A}_5} \to {\cal A}_4$ is dominant. Hence the family $p_i: \textcolor{Red}{\textcolor{Red}{{\cal U}}_i} \to \textcolor{Red}{{\cal A}_5}$ dominates the universal family over \textcolor{Red}{${\cal A}_4$} via the natural map. \endproof \rm \par \noindent
Now let us consider the forgetful map
$$
f_i: \textcolor{Red}{{\cal W}_i} \to \textcolor{Red}{\textcolor{Red}{{\cal U}}_i}.
$$
The fibre of $f_i$ at the moduli point of $(\pi, M)$ is open in  the Prym-Brill-Noether scheme $W^2(\pi)$. Therefore the general Prym-Brill-Noether theory can be applied to such a fibre and to $\textcolor{Red}{{\cal W}_i}$.  In the next proposition we summarize what this theory implies in our particular situation:
\proposition
(0) Each irreducible component of $\textcolor{Red}{{\cal W}_i}$ has dimension $\geq dim \ \textcolor{Red}{\textcolor{Red}{{\cal U}}_i} + 1$. 
 \par \noindent 
(1) $W^2(\pi)$ is non empty and any  of its irreducible components has dimension $\geq 1$.
\par \noindent
(2) The Prym-Petri map $\mu^-_L$ has corank   $\geq 1$ for each $\pi$ and $L \in W^2(\pi)$. \par \noindent
(3) $\mu^-_L$  is injective, equivalently  of corank $1$, \ for $\pi$ general and each $L \in W^2(\pi)$. \par  \noindent 
(4) $W^2(\pi)$ is connected for each $\pi$ and it is a smooth, irreducible curve for a general $\pi$.  
\endproposition \rm \par \noindent
For the proof see [W]. As a standard consequence of  the above properties we have:
\proposition $f_i$ is dominant and $\textcolor{Red}{{\cal W}_i}$ is irreducible, ($i = 1,2$). \endproposition
\proof  Assume $(\pi, M)$ defines a general point $u$ of $\textcolor{Red}{\textcolor{Red}{{\cal U}}_i}$. Then (1) implies that $W^2(\pi)$ is non empty. Since $\pi$ is general, (3) implies that $f_i^*(u) = W^2(\pi)$. Hence $f_i$ is dominant. \\
Let $Y$ be any irreducible component of $\textcolor{Red}{{\cal W}_i}$ and let $y \in Y$ be the moduli point of $(\pi, M, L)$. Since $y \in \textcolor{Red}{{\cal W}_i}$, $\mu^-_L$ is injective i.e. its corank is $1$. This is the dimension of the tangent space at $L$ to the fibre  $f_i^*f_i(y) \subseteq W^2(\pi)$. (1) implies that $f_i^*f_i(y)$ is a smooth curve at $L$ and (0) implies that $dim \ Im (df_y) \geq dim \ \textcolor{Red}{\textcolor{Red}{{\cal U}}_i}$. Since $(\pi, \mu)$ is general it defines a smooth point of \textcolor{Red}{$\textcolor{Red}{{\cal A}_5}$}, hence $dim \  T_{\textcolor{Red}{{\cal W}_i},y}$ $=$ $dim \ \textcolor{Red}{\textcolor{Red}{{\cal U}}_i} + 1$ and $\textcolor{Red}{{\cal W}_i}$ is smooth at $y$. It follows that any irreducible component
$Y$ of $\textcolor{Red}{{\cal W}_i}$ is smooth and dominates $\textcolor{Red}{\textcolor{Red}{{\cal U}}_i}$ via $f / Y$. By (4) a general fibre of $f$ is irreducible, hence $Y = \textcolor{Red}{{\cal W}_i}$. \endproof \rm \par \noindent
For $i = 1,2$ let $x \in \tilde {\mathbb C}_i$,  then the fibre map $\emph q_x:  \textcolor{Red}{\tilde {\cal C}}_x \to \textcolor{Red}{{\cal C}_{x}}$ is an  \'etale double covering. Moreover let $M_x = {\cal O}_{ \textcolor{Red}{\tilde {\cal C}}_x} \otimes \textcolor{Red}{{\cal M}}$ and $L_x = {\cal O}_{ \textcolor{Red}{\tilde {\cal C}}_x} \otimes \textcolor{Red}{{\cal M}}$, by propositions 6.1 and 5.4,  the triple $(\emph q_x, M_x, L_x)$ defines a point of $\textcolor{Red}{{\cal W}_i}$.  Then for $i = 1,2$ the  triple $(\emph q: \textcolor{Red}{\tilde {\cal C}} \to \textcolor{Red}{{\cal C}},\textcolor{Red}{{\cal M}}, \textcolor{Red}{{\cal L}})$ defines a natural morphism
$$
\phi_i: \tilde {\mathbb C}_i \to \textcolor{Red}{{\cal W}_i}
$$  
sending $x \in \tilde {\mathbb C}_i$ to the moduli point of the triple $(q_x, M_x. L_x)$.  Now let us fix a general point
$$
z \in \textcolor{Red}{{\cal W}_i},
$$
our goal will be to show that $ z \in \phi (\tilde {\mathbb C}_i)$. We assume that $z$  is the moduli point of  the triple $$(\pi: \tilde C \to C, M, L). $$  
Let $\mu: H^0(L) \otimes H^0(i^*L) \to H^0 (\omega_{\tilde C})$ be the Petri map of $L$, 
then $Im \ \mu$ defines a map
$$
f_{\mu}: \tilde C \to \mathbf P^2 \times \mathbf P^2 \subset \mathbf P^8.
$$
Up to a projective isomorphism it is not restrictive to assume  $\iota \cdot f_{\mu} = f_{\mu} \cdot i$.
\lemma For a general $z \in \textcolor{Red}{{\cal W}_i} $ the map $f_{\mu}$ is an embedding, ($i = 1,2$). \endlemma
\proof The condition that $f_{\mu}$ is an embedding  is open. Since it holds on $\phi _i (\tilde {\mathbb C}_i)$ it is non empty.
\endproof \rm \par \noindent
From now on we will assume  that $f_{\mu}$ is an embedding and put $\tilde C = f_{\mu}(\tilde C)$ so that
$$
\tilde C \subset \mathbf P^2 \times \mathbf P^2 \subset \mathbf P^8.
$$
\lemma For a general $z \in \textcolor{Red}{{\cal W}_i}$ and $i = 1,2$ we have: \par \noindent
(1) there exists a unique hyperplane containing $\tilde C$ and its equation is $\iota^*$-invariant. \par \noindent
(2) Such a hyperplane is transversal to $\mathbf P^2 \times \mathbf P^2$ and to its diagonal $\Delta$.
\endlemma
\proof (1) $Ker \ \mu$ is naturally isomorphic to the space of linear forms vanishing on $\tilde C$, let us show that $dim \ Ker \ \mu \geq 1$ for each $z \in \textcolor{Red}{{\cal U}}$. This follows from Brill-Noether theory, indeed note that
$$
L \in W^2(\pi) \subseteq W := \lbrace N \in Pic^8 ( \tilde C) \ / \ h^0(N) \geq 3 \rbrace.
$$
$W$ is a Brill-Noether locus and the tangent space at $L$  to $W$ is $Coker \ \mu$. On the other hand we know that $dim \ W^2(\pi) \geq 1$, hence $dim \ Coker \ \mu \geq 1$ for each $z \in \textcolor{Red}{{\cal U}}$. Since $\mu$ is a linear map between vector spaces of the same dimension, $dim \ Ker \ \mu \geq 1$. We know that $dim \ Ker \ \mu = 1$ if $z \in \phi_i(\tilde {\mathbb C})$. By semicontinuity the same is true for a general $z$.  Finally $i^*$ acts as an involution on $Ker \ \mu$. Since the Prym-Petri map is injective no $-1$ vector is in $Ker \ \mu$ hence $h = i^*h$,  $\forall h \in Ker \ \mu$. \\
(2) The transversality condition  is open and non empty  because each $z \in \phi_i( \tilde {\mathbb C}_i)$ satisfies it.   \endproof \rm \par \indent
Cutting $\mathbf P^2 \times \mathbf P^2$ with the hyperplane $< \tilde C>$ we obtain a smooth threefold
$$
\tilde T \ = \ < \tilde C >  \ \cap \ \mathbf P^2 \times \mathbf P^2.
$$
\lemma For $i = 1,2$ and a general $z \in \textcolor{Red}{{\cal W}_i}$ one has:
\par \noindent
(1) $ h^0 ({\cal I}_{\tilde C / \tilde T}(2)) = 3$ and the equation of each $\tilde S \in \mid {\cal I}_{\tilde C / \tilde T} (2)\mid$ is $\iota^*$-invariant.
\par \noindent
(2) A general $\tilde S \in \mid {\cal I_{\tilde C / \tilde T}}(2) \mid$ is smooth and transversal to the diagonal $\Delta$.
\par \noindent
(3) A general  $\tilde S$ contains a smooth conic $B$ such that $B \cap \iota^*B = \emptyset$, moreover 
$$ \mid \tilde C - \tilde F \mid \ = \  \mid  \tilde B +  \iota^*B \mid, $$ where $\tilde F$ is a hyperplane section of $\tilde S$.
\endlemma
\proof (1) By 5.3  ${\cal I}_{\tilde C / \tilde T}$ is acyclic and $h^0({\cal I}_{\tilde C / \tilde T}(2))= 3$ if $z \in \phi_i (\tilde {\mathbb C}_i)$. By semicontinuity the same is true for a general $z$. $\iota^*$ acts as an involution on $H^0({\cal I}_{\tilde C / \tilde T}(2))$: by 5.3 its $-1$ eigenspace is zero if  $z \in  \phi( \tilde {\mathbb C})$. So this is generically true on $\textcolor{Red}{{\cal U}}$ and hence $\iota^*$ is the identity for a general $z$.   \\ 
(2) If $z \in \phi_i (\tilde {\mathbb C}_i)$ we know that $\mid {\cal I}_{\tilde C / \tilde T}(2)\mid$ has minimal dimension two and that
its general element $\tilde S$ is smooth and transversal to $\Delta$. Then the same is true for a general $z \in \textcolor{Red}{{\cal U}}$. \par \noindent
(3) Let $\tilde S \in \mid {\cal I}(2) \mid$ be general and let $\tilde F$ be a hyperplane section of $\tilde S$:  at first we want to show  that $\mid \tilde C - \tilde F\mid \neq \emptyset $. For this observe that ${\cal O}_{\tilde C}(\tilde F) \cong \omega_{\tilde C}$ and consider the standard exact sequence
$$
0 \to {\cal O}_{\tilde S}(\tilde F - \tilde C) \to {\cal O}_{\tilde S}(\tilde F) \to {\cal O}_{\tilde C}(\tilde F) \to 0.
$$
Since $h^0({\cal O}_{\tilde S} (\tilde F - \tilde C)) = h^1({\cal O}_{\tilde S}(\tilde F))= 0$  the associated long exact sequence is
$$
0 \to H^0({\cal O}_{\tilde S}(\tilde F)) \to H^0(\omega_{\tilde C}) \to H^1({\cal O}_{\tilde S}(\tilde F - \tilde C)) \to 0,
$$
hence  $h^1({\cal O}_{\tilde S}(\tilde F - \tilde C)) = 1$. From Riemann-Roch and Serre duality  it follows $h^2({\cal O}_{\tilde S}(\tilde F - \tilde C))$ $ = h^0({\cal O}_{\tilde S}(\tilde C - \tilde F)) = 1$ and $\mid \tilde C - \tilde F \mid \neq \emptyset$. \   If $z \in \phi (\tilde {\mathbb C})$ then  $\mid \tilde C - \tilde F \mid$ consists of one element $B + \iota^*B$, where $B$ is a smooth conic and $B \cap \iota^*B = \emptyset$. So the same is true for a general $z$. \endproof
\rm \par \noindent
Since $z$ is the moduli point of the triple $(\pi, \tilde C \to C, M, L)$ we now study $M$. At first let us recall that $h^0(M) = i$, since $z \in \textcolor{Red}{{\cal W}_i}$. Moreover $M = {\cal O}_{\tilde C}(d)$ where $d$ is an effective divisor of degree 8 such that $\pi_*d \in \mid \omega_C \mid$,  of  course we have
$$
d \subset \tilde C \subset  \tilde S \subset \mathbf P^2 \times \mathbf P^2.
$$
\lemma Let $z \in \textcolor{Red}{{\cal W}_i}$ be general, let $i = 1,2$ and let  ${\cal }I_{d / \tilde S}$ be the ideal sheaf of $d$ in $\tilde S$, then: \par \noindent
(1) $d$ is contained in a unique hyperplane section $\tilde F$ of $\tilde S$. \par \noindent
(2) $\tilde F$ is smooth, transversal to $B + \iota^*B$ and its equation is $\iota^*$-invariant.
\par \noindent 
(3) $\mid d \mid$ is a base-point-free pencil on $\tilde F$. \par \noindent
(4) $d \sim a + b' + b''$ where $b', b'' \in \tilde F \cap \ (B \cup \ \iota^*B)$ and
$b'' \neq \iota(b')$. \par \noindent
(5) $z \in \textcolor{Red}{{\cal W}_1}$  if $\# (\lbrace b' + b'' \rbrace  \cap B)$ is odd and $z \in \textcolor{Red}{{\cal W}_2}$ if  $\# (\lbrace b' + b'' \rbrace  \cap B)$ is even. \par \noindent
(6) the linear system $P_z := \mid {\cal I}_{d / \tilde S}(\tilde C) \mid^+$ is a pencil of smooth, irreducible curves. \par \noindent
If $o = (y, B)$, where $y = (a,b',b'')$,  it follows that $(o, \tilde S, d)$ is a marking 0-cycle of type $i$. 
\endlemma
\proof (1) Note that, by 6.7 (1) the space $<\tilde S>$ is $\iota$-invariant, in particular the projectivized eigenspaces of $\iota / < \tilde S >$ are $\mathbf P^-$ and $\mathbf P^{4+}$ $:=$ $ < \tilde S >  \cap \ \mathbf P^+$. Notice also that $\tilde C$ is not in a hyperplane of $< \tilde S >$ for degree reasons.
On the other hand we know that $q / \tilde C = \pi$ and $\pi_* d \in \mid \omega_C \mid$. Hence it follows that $q(\tilde C)$ is $C$  canonically embedded in $\mathbf P^{4+}$ and that $<q_*d>$ is a hyperplane in $\mathbf P^{4+}$. This implies that $d$ is always contained in the hyperplane
$$
(q / < \tilde S >)^* < q_*d >,
$$
whose equation is $\iota^*$-invariant. This implies that $h^0({\cal I}_{d / \tilde S}(1))^+ \geq 1$ for
each $z \in \textcolor{Red}{{\cal W}_i}$. By Proposition 4.3 the equality holds if $z \in \phi_i(\tilde {\mathbb C}_i))$. Hence, by semicontinuity, it holds generically on $\textcolor{Red}{{\cal W}_i}$. \\
(2) We have shown in (1) that $< d >$ is a hyperplane whose equation is $\iota^*$-invariant. Actually $< d >$ is the hyperplane of $< \tilde S>$ spanned by $\mathbf P^-$ and $< q_*d>$. Since $z$ is general we can assume that $d$ is general in the irreducible family $D_i$ of the effective divisors  $d'$ on $ \tilde C$  such that $\pi_*d' \in \mid \omega_C \mid$ and $h^0({\cal O}_{\tilde C}(d')) = i$. It is well known that the image of $D_i$ via the push-down map $\pi_*$ is open in $\mid \omega_C \mid$, (cfr. [B1]).
 Therefore, for a general $d' \in P'$, $< d' >$ is just a general hyperplane through $\mathbf P^-$. Hence it is transversal to $\tilde S$ and $B + \iota^*B$ and the same holds for $< d >$. \\
(3) Obviously no hyperplane of $< d >$ contains $d$.  Since $\tilde F$ is canonically embedded in $< d >$ it follows that $d$ is a non special divisor and that $h^0({\cal O}_{\tilde F}(d)) = 2$. To see that $\mid d \mid$ is base-point-free observe that this is an open condition, which is  satisfied on $\phi_i ( \tilde {\mathbb C}_i)$ by Proposition 4. 3.\\
(4) As usual we put $C = q(\tilde C)$, $F = q(\tilde F)$, $\overline B = q(B + \iota^*B)$ and $S = q(\tilde S)$. Then we consider the double covering $q/\tilde S: \tilde S \to S$ and observe that  $ C \sim F + \overline B$. On the other hand $\pi_*d = q_*d$ is a hyperplane section of $C$, so we can conclude that $\pi_*d = F \cdot C$. Since $C \sim F + \overline B$, it follows that
$$
\pi_*d \in  \mid \omega_F (\overline B) \mid = \mid \omega_F(\pi(b_1) + \pi(b_2)) \mid.
$$
Since $F$ is non hyperelliptic the divisor $\pi(b_1) + \pi(b_2)$ is isolated on $F$, it is also smooth because $\tilde F$ is transversal to $B + \iota^*B$ so that $F$ is transversal to $\overline B$. We consider the curve
$$
\Gamma = \lbrace \pi_*m, \ m \in \mid {\cal O}_{\tilde F}(d) \mid \rbrace \subset \mid \omega_F(\overline B) \mid.
$$
As is well known the divisors of  $\mid \omega_F(\pi(b_1) + \pi(b_2)) \mid$ passing through $\pi(b_1) + \pi(b_2)$ form
a hyperplane $H$. Hence $H \cap \Gamma$ is non empty and this implies that $d \sim \pi_* (a + b' + b'')$ for some $a \in Div \  \tilde F$ such that $q_*a \in \mid \omega_F \mid$
and for some $b', b'' \in B \cup \iota^*B$ such that $b'' \neq \iota(b')$.
Since $\mid {\cal O}_{\tilde F}(d) \mid$ is base-point-free $h^0({\cal O}_{\tilde F}(a)) = 1$, hence the linear span of $a$ is a 5-space $\mathbf P^5_a$. Notice also that $\mathbf P^- \subset \mathbf P^5_a$ because $q_*a$ is a plane. Hence it follows that $q_*a \in \mathbb A$ for a general $z$. \\ Note that, by the previous part of the proof,  assumptions (1), (2), (3) of parity lemma are satisfied, while (4) is obvious because $d$ is contained in the smooth curve $\tilde C$.  Then properties (5) and (6) follow from parity lemma 5.4 and proposition 5.5 (1).
\endproof \rm \par \noindent
Putting together the previous lemmas we can finally deduce that:  
\lemma The map $\phi_i: \tilde {\mathbb C}_i \to \textcolor{Red}{{\cal W}_i}$ is dominant for $i = 1,2$. \endlemma
\proof Let $z$ be general in $\textcolor{Red}{{\cal W}_i}$, then $z$ is the moduli point of $(\pi: \tilde C \to C, M, L)$,
where $M ={\cal O}_{\tilde C}(d)$ and $d$ is a smooth, effective divisor of degree 8 as above. Due to the \vfill \eject \noindent previous lemmas we have the  embeddings
$$
d \subset \tilde C \subset  \tilde S \subset  \mathbf P^2 \times \mathbf P^2 \subset \mathbf P^8,
$$
where $L \cong {\cal O}_{\tilde C}(1,0)$ and $i = \iota / \tilde C$. Moreover we have the smooth conic 
$B \subset \tilde S$ such that $\tilde C - \tilde F \sim B + \iota^*B$ and the smooth hyperplane section 
$\tilde F = < d > \cap \  \tilde S$. We have also shown in the above Lemma 6.8 that, putting $y = (a, b', b'')$ and  $o = (y, B)$,  it follows that $(o, \tilde S, d)$ is a marking $0$-cycle of type $i$, hence
$
x = (o, \tilde S, d, \tilde C) \in \tilde {\mathbb C}_i.
$
On the other hand $\phi_i(x)$ is the moduli point of the triple $(q_x, M_x, L_x)$, where $q_x: \tilde C \to \tilde C / < \iota>$ is the quotient map, $M_x = {\cal O}_{\tilde C}(d)$ and $L_x = {\cal O}_{\tilde C}(1,0)$. Therefore $\phi_i(x) = z$ and $\phi_i$ is dominant. \endlemma
\theorem The universal principally polarized abelian variety over \textcolor{Red}{${\cal A}_4$} and the universal theta divisor over \textcolor{Red}{${\cal A}_4$} are unirational varieties.
\endtheorem
\proof We have seen that $\textcolor{Red}{{\cal W}_1}$ dominates the universal p.p.a.v. over \textcolor{Red}{${\cal A}_4$} and that $\textcolor{Red}{{\cal W}_2}$ dominates the universal theta divisor over \textcolor{Red}{${\cal A}_4$}. On the other hand $\phi_i:  \tilde {\mathbb C}_i \to \textcolor{Red}{{\cal W}_i}$ is
dominant and $\tilde {\mathbb C}_i$ is rational. Hence the universal p.p.a.v. and the universal theta divisor over \textcolor{Red}{${\cal A}_4$} are unirational varieties.  \endproof
 \bibliographystyle{amsalpha}

 \rm \par \noindent
 ADDRESS Alessandro Verra, Universit\'a Roma Tre, Dipartimento di Matematica \\
Largo San Leonardo Murialdo 1, \ 00146 Rome, Italy \\
 E-MAIL: \ verra@mat.uniroma3.it
    \enddocument